\documentclass[12pt]{amsart}

\usepackage{amssymb}
\usepackage{amscd}
\usepackage{statex}

\title{Galois lines for the Giulietti--Korchm\'{a}ros curve}
\author{Satoru Fukasawa \and Kazuki Higashine}

\subjclass[2010]{14H50, 14H37}
\keywords{Galois point, Galois line, Giulietti--Korchm\'{a}ros curve, automorphism group}
\thanks{The first author was partially supported by JSPS KAKENHI Grant Number 16K05088}
\address[Satoru Fukasawa]{Department of Mathematical Sciences, Faculty of Science, Yamagata University, Kojirakawa-machi 1-4-12, Yamagata 990-8560, Japan}
\email{s.fukasawa@sci.kj.yamagata-u.ac.jp}
\address[Kazuki Higashine]{Graduate School of Science and Engineering, Yamagata University, Kojirakawa-machi 1-4-12, Yamagata 990-8560, Japan}
\email{s16m111m@st.yamagata-u.ac.jp}

\newtheorem{thm}{Theorem}

\newtheorem{lem}{Lemma}

\newtheorem{fact}{Fact}

\theoremstyle{definition}
\newtheorem{rem}{Remark}

\begin{document}
\begin{abstract}
We describe the arrangement of all Galois lines for the Giulietti--Korchm\'{a}ros curve in the projective $3$-space. 
As an application, we determine the set of all Galois points for a plane model of the GK curve. 
This curve possesses many Galois points.  
\end{abstract}
\maketitle

\section{Introduction} 
In Algebraic Geometry, the Hermitian curve 
$$ \mathcal{H}: x^q+x-y^{q+1}=0$$
in the projective plane $\mathbb{P}^2$ over a field of characteristic $p>0$ has many interesting and important properties, where $q$ is a power of $p$. 
The following beautiful theorem presented by Homma in the theory of Galois point represents one of them (see \cite{miura-yoshihara, yoshihara, yoshihara-fukasawa} for the definition of the Galois point): 

\begin{fact}[\cite{homma}]
For the Hermitian curve $\mathcal{H}$, the set of all Galois points coincides with the set of all $\mathbb{F}_{q^2}$-rational points of $\mathbb{P}^2$. 
\end{fact}
 
It would be good to obtain a result similar to the theorem of Homma. 
In this article, we focus on the Giulietti--Korchm\'{a}ros curve. 

The Giulietti--Korchm\'{a}ros curve $\mathcal{X}$ is an important class of maximal curves over a finite field, which is described as a complete intersection of two hypersurfaces 
$$ x^q+x-y^{q+1}=0, \mbox{ and } y((x^q+x)^{q-1}-1)-z^{q^2-q+1}=0 $$
in $\mathbb{P}^3$ (\cite{giulietti-korchmaros}).  
Our main result is the following theorem on the arrangement of ``Galois lines'' in $\mathbb{P}^3$ for the GK curve $\mathcal{X}$ (i.e. the projection from the line induces a Galois extension, see \cite{duyaguit-yoshihara, yoshihara2}). 

\begin{thm} \label{galois_line} 
The set of all Galois lines for $\mathcal{X}$ coincides with the set of all $\mathbb{F}_{q^2}$-lines $\ell$ with $\ell \ni (0:0:1:0)$ or $\ell \subset \{Z=0\}$.  
\end{thm}

Giulietti and Korchm\'{a}ros introduced a plane model  
$$ x^{q^3}+x-(x^q+x)^{q^2-q+1}-z^{q^3+1}=0 $$
(\cite[Theorem 4]{giulietti-korchmaros}). 
The projective closure $\mathcal{X}'$ of this curve is the same as the image $\pi_R(\mathcal{X})$ under the projection $\pi_R:\mathbb{P}^3 \dashrightarrow \mathbb{P}^2$ from the point $R:=(0:1:0:0)$.   
As an application of Theorem \ref{galois_line}, we will obtain the following. 

\begin{thm}\label{galois_point} 
The set of all Galois points for $\mathcal{X}'$ coincides with the set $$\{(0:1:0)\} \cup (\mathcal{X}' \cap \{Z=0\}).$$ 
\end{thm}

It is important to find a new example of a plane curve with many Galois points, in the theory of Galois point (see \cite{yoshihara-fukasawa}). 
According to Theorem \ref{galois_point}, the number of Galois points in $\mathcal{X}' \setminus {\rm Sing}(\mathcal{X}')$ is exactly $q+1$, where ${\rm Sing}(\mathcal{X}')$ is the set of all singular points of $\mathcal{X}'$.

\section{Galois lines with degree $q^3$} 

We consider the Giulietti--Korchm\'{a}ros curve over an algebraically closed field $k$ of characteristic $p>0$. 
We take systems $(X:Y:Z:W)$ and $(X:Z:W)$ of homogeneous coordinates of $\mathbb{P}^3$ and of $\mathbb{P}^2$ respectively. 
Let $x=X/W$, $y=Y/W$ and $z=Z/W$. 
For points $P$ and $Q \in \mathbb{P}^3$ with $P \ne Q$, the line passing through $P$ and $Q$ is denoted by $\overline{PQ}$.  

Let $P_{\infty}:=(1:0:0:0) \in \mathcal{X}$ and let $\ell_{\infty} \subset \mathbb{P}^3$ be the line defined by $Z=W=0$. 
Then, $P_{\infty} \in \ell_{\infty}$, and $\ell_{\infty}=\overline{RP_{\infty}}$. 
First, we show that $\ell_{\infty}$ is a Galois line.  
Note that the affine part $\mathcal{X}_W$ of $\mathcal{X}$ given by $W \ne 0$ is the same as the curve defined by 
$$ x^q+x-y^{q+1}=y^{q^2}-y-z^{q^2-q+1}=0. $$ 
The subgroup 
$$G_1:=\left\{
\left(\begin{array}{cccc}
1 & b^q &  0 &a \\
0 & 1 &  0 &b \\
0 & 0 & 1 & 0 \\
0 & 0 & 0 & 1
\end{array}\right)
\ | \ a, b \in \mathbb{F}_{q^2}, \ a^q+a-b^{q+1}=0\right\}
\subset {\rm PGL}(4, k)$$ 
of order $q^3$ acts on $\mathcal{X}$ (\cite[p. 238]{giulietti-korchmaros}). 
It is not difficult to check that $k(\mathcal{X})^{G_1}=k(z)$. 
The extension $k(\mathcal{X})/k(\mathcal{X})^{G_1}$ coincides with the extension $k(\mathcal{X})/k(z)$ induced by the projection $\pi_{\ell_{\infty}}: \mathbb{P}^3 \dashrightarrow \mathbb{P}^1; (X:Y:Z:W) \mapsto (Z:W)$ from $\ell_{\infty}$. 
Therefore, $\ell_{\infty}$ is a Galois line with $\ell_{\infty} \cap \mathcal{X}=\{P_{\infty}\}$.  

Note that $\ell_{\infty}$ coincides with the tangent line of the Hermitian curve $Z=X^qW+XW^q-Y^{q+1}=0$ at $P_{\infty}$.  
We infer that there exist $q^3+1$ Galois lines for $\mathcal{X}$, by the following theorem of Giulietti and Korchm\'{a}ros. 

\begin{fact}[A weak version of Theorem 7 in \cite{giulietti-korchmaros}] \label{action} 
The full automorphism group ${\rm Aut}(\mathcal{X})$ of $\mathcal{X}$ acts on $\mathcal{X} \cap \{Z=0\}$, and the action of the subgroup ${\rm Aut}(\mathcal{X}) \cap {\rm PGL}(4, k)$ on this set is doubly transitive. 
\end{fact}

We would like to show that the number of lines which induce a Galois extension of degree $q^3$ is at most $q^3+1$. 
The following two lemmas are needed. 

\begin{lem}[see \cite{beelen-montanucci}, p.13] \label{order} 
Let $P \in \mathcal{X}$ and let $H \subset \mathbb{P}^3$ be a hyperplane with $H \ni P$.  
\begin{itemize}
\item[(1)] If $P \in \mathcal{X} \cap \{Z=0\}$, then ${\rm ord}_{P}H= 1$, $q^2-q+1$ or $q^3+1$. 
\item[(2)] If $P \not\in \mathcal{X} \cap \{Z=0\}$, then ${\rm ord}_{P}H=1$, $q$, $q^3$ or $q^3+1$. 
\end{itemize}
\end{lem} 

\begin{lem} \label{hermitian} 
Let $\mathcal{H}$ be the Hermitian curve given by $Z=XW^q+XW^q-Y^{q+1}=0$ and let $\mathcal{H}(\mathbb{F}_{q^2})$ be the set of all $\mathbb{F}_{q^2}$-rational points on $\mathcal{H}$. 
For each line $\ell \subset \{Z=0\}$, $\#(\mathcal{H}(\mathbb{F}_{q^2}) \cap \ell)=0, 1$ or $q+1$. 
\end{lem}

Assume that a line $\ell \subset \mathbb{P}^3$ induces a Galois extension of degree $q^3$. 
Note that the Galois group $G_{\ell}$ of order $q^3$ acts on the set $\mathcal{X} \cap \{Z=0\}$ of cardinality $q^3+1$, by Fact \ref{action}. 
By a fact of group theory (see \cite[Chapter 2, Section 1 (1.3)]{suzuki}), there exists a point $P \in \mathcal{X} \cap \{Z=0\}$ fixed by any element of $G_{\ell}$.    
It follows from \cite[III. 8.2]{stichtenoth} that the ramification index at $P$ is equal to $q^3$ for the projection $\pi_{\ell}$ from $\ell$.  
We can assume that $P=P_{\infty}$. 
If $P_{\infty} \not\in \ell$, then there exists a hyperplane $H$ such that ${\rm ord}_{P_{\infty}}H=q^3$. 
This is a contradiction to Lemma \ref{order}.
Therefore, $P_{\infty} \in \ell$. 
Since $\pi_{\ell}^{-1}(\pi_{\ell}(P_{\infty}))=\{P_{\infty}\}$, the hyperplane $W=0$ includes $\ell$. 
Assume that $\ell \not\subset \{Z=0\}$. 
Then, $\ell \cap \{Z=0\}=\{P_{\infty}\}$ and $\#(\mathcal{X} \cap \{Z=0\} \cap H)=1$ or $q+1$ for each hyperplane $H \supset \ell$, by Lemma \ref{hermitian}.
By Fact \ref{action} and \cite[III. 7.2]{stichtenoth}, there exists a ramification point different from $P_{\infty}$ with index a power of $p$. 
This is a contradiction to Lemma \ref{order}. 
Therefore, $\ell$ is defined by $Z=W=0$. 
The proof of the assertion for Galois lines with $[k(\mathcal{X}):k(\mathbb{P}^1)]=q^3$ in Theorem \ref{galois_line} is completed.

We consider Galois points in $\mathcal{X}' \setminus {\rm Sing}(\mathcal{X}')$. 
Note that Galois lines $\ell$ for $\mathcal{X}$ passing through $R=(0:1:0:0)$ correspond to Galois points $P \in \mathbb{P}^2$ for $\mathcal{X}'$, by $\ell \mapsto \pi_R(\ell)$, since the projection $\pi_R: \mathcal{X} \rightarrow \mathcal{X}'$ is birational.  
Since all Galois lines $\ell$ with degree $q^3$ are included in the plane $\{Z=0\}$, all Galois points with degree $q^3$ on $\mathcal{X}'$ are contained in the line $\{Z=0\}$ in $\mathbb{P}^2$.  
Note that if a line $\ell$ with $R \in \mathcal{\ell} \subset \{Z=0\}$ intersects $\mathcal{X}$ at $q+1$ points, then $\pi_R(\ell) \in {\rm Sing}(\mathcal{X})$.    
Considering $\mathbb{F}_{q^2}$-lines passing through $R$, we infer that there exist exactly $q+1$ Galois points on $\mathcal{X}'$ with degree $q^3$, which points are $(\alpha:0:\beta) \in \mathbb{P}^2$ with $\alpha^q\beta+\alpha\beta^q=0$. 
The assertion in Theorem \ref{galois_point} for Galois points with degree $q^3$ follows. 

\section{Galois lines with degree $q^3+1$}

We consider Galois lines $\ell \subset \mathbb{P}^3$ with $\ell \cap \mathcal{X} =\emptyset$. 
Let $R'=(0:0:1:0)$ and let $\ell_0$ be the line defined by $X=W=0$. 
Then, $\ell_0 \cap \mathcal{X}=\emptyset$, and $\ell_{0}$ coincides with the line passing through $R'$ and $R=(0:1:0:0)$. 
The subgroup 
$$G_2:=\left\{
\left(\begin{array}{cccc}
1 & 0 &  0 &0 \\
0 & \zeta^{q^2-q+1} &  0 &0 \\
0 & 0 & \zeta & 0 \\
0 & 0 & 0 & 1
\end{array}\right)
\ | \ \zeta \in k, \ \zeta^{q^3+1}=1\right\}
\subset {\rm PGL}(4, k)$$ 
of order $q^3+1$ acts on $\mathcal{X}$. 
It is not difficult to check that $k(\mathcal{X})^{G_2}=k(x)$. 
The extension $k(\mathcal{X})/k(\mathcal{X})^{G_2}$ coincides with the extension $k(\mathcal{X})/k(x)$ induced by the projection $\pi_{\ell_0}: \mathbb{P}^3 \dashrightarrow \mathbb{P}^1; (X:Y:Z:W) \mapsto (X:W)$ from $\ell_{0}$. 
Therefore, $\ell_{0}$ is a Galois line with $\ell_{0} \cap \mathcal{X}=\emptyset$.  

By Fact \ref{action}, for each $\mathbb{F}_{q^2}$-rational point $Q \in \{Z=0\} \setminus \mathcal{X} \subset \mathbb{P}^3$, there exists a line $\ell \subset \mathbb{P}^3$ with $ Q \in \ell \not\subset \{Z=0\}$ such that $\pi_{\ell}$ induces a Galois extension of degree $q^3+1$. 
Therefore, the number of lines with required properties is at least $q^4+q^2+1-(q^3+1)=q^4-q^3+q^2$.   

We consider the case where a Galois line $\ell$ is included in the plane $Z=0$. 
Note that the projection $\pi_{\ell}$ is not ramified at each points in $\mathcal{X} \cap \{Z=0\}$. 
By Lemma \ref{order} and \cite[III. 7.2]{stichtenoth}, the ramification index at all ramification points for $\pi_{\ell}$ is equal to $q^3+1$. 
By the Riemann--Hurwitz formula, the integer $2g-2+2(q^3+1)$ is divisible by $q^3$, for the genus $g$ of $\mathcal{X}$. 
This is a contradiction to \cite[Theorem 2]{giulietti-korchmaros} (this integer is equal to $(q^3+1)q^2$). 
Therefore, $\ell \not\subset \{Z=0\}$ holds for all Galois lines $\ell$ with $\ell \cap \mathcal{X}=\emptyset$. 

Let $\ell \subset \mathbb{P}^3$ be a line with $\ell \cap \{Z=0\}=\{Q\} \not\subset \mathcal{X}$ which induces a Galois extension of degree $q^3+1$. 
It follows from Lemma \ref{hermitian} that $\#(\mathcal{X} \cap \{Z=0\} \cap  H)=0$, $1$ or $q+1$, for each hyperplane $H \supset \ell$.     
By Fact \ref{action} and \cite[III. 7.2]{stichtenoth}, $\pi_{\ell}$ is ramified at points in $\mathcal{X} \cap \{Z=0\}$ with index $q^3+1$ or $q^2-q+1$. 
Note that if the index at $P$ is $q^3+1$, then the line $\overline{QP}$ is $\mathbb{F}_{q^2}$-rational. 
Considering lines in the plane $Z=0$ passing through $Q$, there exist at least two lines  over $\mathbb{F}_{q^2}$ containing $Q$. 
Therefore, the point $Q$ is an $\mathbb{F}_{q^2}$-rational point.  
By Fact \ref{action}, we can assume that $Q=R=(0:1:0:0)$. 

We would like to show the uniqueness of the Galois line $\ell \ni R$ with $\ell \cap \mathcal{X}=\emptyset$. 
We consider the projection $\pi_R: \mathbb{P}^3 \dashrightarrow \mathbb{P}^2$. 
Since all points of $\mathcal{X} \cap \{Z=0\}$ are ramification points for $\pi_\ell$, all tangent lines at smooth points in $\mathcal{X}' \cap \{Z=0\}$ contains the point $\pi_R(\ell)$. 
This implies that $\pi_R(\ell)=(0:1:0)$. 
The uniqueness follows. 
This observation also implies that the number of outer Galois points for $\mathcal{X}'$ is exactly one.

\begin{rem}
The tangent line at each point of $\mathcal{X} \cap \{Z=0\}$ passes through $R'=(0:0:1:0)$ (see \cite[p.234]{giulietti-korchmaros}). 
Furthermore, $R' \in \ell$ for all Galois lines $\ell$ with $\ell \cap \mathcal{X} =\emptyset$. 
Therefore, a Galois line $\ell$ with $\ell \cap \mathcal{X}=\emptyset$ coincides with the line passing through $R'$ and an $\mathbb{F}_{q^2}$-rational point in the plane $Z=0$ but not on $\mathcal{X}$. 
\end{rem}

\begin{rem}
For each Galois line $\ell$ with $\ell \cap \mathcal{X}= \emptyset$, the Galois group $G_\ell$ includes the subgroup 
$$ \left\{
\left(\begin{array}{cccc}
1 & 0 &  0 &0 \\
0 & 1 &  0 &0 \\
0 & 0 & \eta & 0 \\
0 & 0 & 0 & 1
\end{array}\right)
\ | \ \eta \in k, \ \eta^{q^2-q+1}=1 \right\}
\subset {\rm PGL}(4, k). $$
Therefore, $G_{\ell} \cap G_{\ell'} \ne \{1\}$ for each two Galois lines $\ell$ and $\ell'$ not intersecting $\mathcal{X}$. 
According to \cite[Theorem 1]{fukasawa}, there exist no plane model of $\mathcal{X}$ realizing $G_{\ell}$ and $G_{\ell'}$ as Galois groups at two outer Galois points. 
\end{rem}

\section{Galois lines with degree at most $q^3-1$}
The tangent line at each point of $\mathcal{X} \cap \{Z=0\}$ is a Galois line. 
In fact, the projection from the line $Y=W=0$ induces the extension $k(x, y, z)/k(y)$, and $y$ is fixed by automorphisms $(x, y, z) \mapsto (x+\alpha, y, \eta z)$, where $\alpha^q+\alpha=0$ and $\eta^{q^2-q+1}=1$.

We show that any line $\ell \subset \{Z=0\}$ such that $\ell \cap \mathcal{X}$ contains at least two points is a Galois line. 
It follows from Lemma \ref{hermitian} that $\ell$ is $\mathbb{F}_{q^2}$-rational and contains exactly $q+1$ points of $\mathcal{X}$. 
We consider the line $Y=Z=0$. 
Then, the extension is $k(x,y,z)/k(y/z)$. 
The automorphisms $\sigma_{\alpha}: (x, y, z) \mapsto (x+\alpha, y, z)$ and $\tau: (x,y, z) \mapsto (\xi^{q+1}x, \xi y, \xi z)$ act on $\mathcal{X}$, where $\alpha^q+\alpha=0$ and $\xi$ is a primitive $(q-1)$-th root of unity, and fix $y/z$. 
Note that the group generated by such $q(q-1)$ automorphisms fixes $P_{\infty}$. 
We consider the linear transformation $\varphi$ on $\mathbb{P}^4$ represented by 
$$ A=\left(\begin{array}{cccc} 
1 & 0 & 0 & \rho^q \\
0 & 1 & 0 & 0 \\
0 & 0 & -1 & 0 \\
1 & 0 & 0 &-\rho 
\end{array}\right), $$
where $\rho^q+\rho=1$ (see \cite{gqz}). 
Then, $\varphi(\mathcal{X})$ is given by 
$$ x^{q+1}-1=y^{q+1}, \mbox{ and } \ y\frac{x^{q^2}-x}{x^{q+1}-1}=z^{q^2-q+1}. $$
The linear transformation $\psi$ given by 
$$ \left(\begin{array}{cccc} \beta^2 & 0 & 0 &0 \\
0 & \beta & 0 & 0 \\
0 & 0 & \beta & 0 \\
0 & 0 & 0 & 1 
\end{array}\right), $$
where $\beta$ is a primitive $(q+1)$-th root of unity, 
acts on $\varphi(\mathcal{X})$. 
Since $\varphi^*(y/z)=-y/z$ and $\psi^*(y/z)=y/z$, $\mu^*:=(\varphi^{-1}\circ\psi\circ\varphi)^*$ fixes $y/z$. 
Let $G_3 \subset {\rm Aut}(\mathcal{X})$ be the group generated by $\tau, \mu$ and all $\sigma_{\alpha}$. 
Since $\beta^2 \ne 1$, $\mu(P_{\infty})=(\beta^2\rho+\rho^q:0:0:\beta^2-1) \ne P_{\infty}$.  
It follows that $G_3$ acts on $\mathcal{X} \cap \{Y=Z=0\}$ transitively. 
Therefore, there exist at least $q(q-1)$ elements of $G_3$ fixing $P$, for each point $P \in \mathcal{X} \cap\{Y=Z=0\}$.  
It follows that $|G_3| \ge q(q-1)(q+1)$ and hence, the line $Y=Z=0$ is a Galois line whose Galois group is equal to $G_3$. 
Furthermore, the image of this line under the projection from an $\mathbb{F}_{q^2}$-rational point in $\{Y=Z=0\} \setminus \mathcal{X}$ is a singular point of the image of $\mathcal{X}$, where the image of $\mathcal{X}$ is projectively equivalent to $\mathcal{X}'$, and this point is also a Galois point. 

\begin{rem}
The automorphism $\varphi^{-1}\circ\psi\circ\varphi$ is represented by 
$$\left(\begin{array}{cccc} 
\beta^2\rho+\rho^q & 0 & 0 & (\beta^2-1)\rho^{q+1} \\
0 & \beta & 0 & 0 \\
0 & 0 & \beta & 0 \\
\beta^2-1& 0 & 0 & \beta^2\rho^q+\rho 
\end{array}\right). $$
\end{rem}

Assume that a line $\ell \subset \mathbb{P}^3$ induces a Galois extension of degree $d \le q^3-1$. 
By Lemma \ref{hermitian} and the previous paragraph, we can assume that $\ell \not\subset \{Z=0\}$.

We consider the case where $\ell \cap \{Z=0\}=\{Q\} \subset \mathcal{X}$. 
We can assume that $Q=P_{\infty}$ and $\ell$ is not the tangent line at $P_{\infty}$. 
Note that $\pi_{P_{\infty}}(\mathcal{X}) \subset \mathbb{P}^2$ is defined by $y^{q^2}-y=z^{q^2-q+1}$, and the projection from each point of $\mathbb{P}^2$ is not birational onto $\mathbb{P}^1$ for this curve.  
Since $\pi_{\ell}$ factors through the projection $\pi_{P_{\infty}}$, it follows that $d >q$. 
By Lemma \ref{hermitian}, for each hyperplane $H \supset \ell$, $\# ((\mathcal{X} \cap\{Z=0\} \cap H)\setminus \{P_{\infty}\})=0$ or $q$. 
Then, there exists a ramification point different from $P_{\infty}$ with index $d/q$. 
It follows from Lemma \ref{order} that $d=q(q^2-q+1)$ or $q$. 
Therefore, $d=q(q^2-q+1)$. 
Let $P_1$ and $P_2 \in \mathcal{X} \cap \{Z=0\} \setminus \{P_{\infty}\}$ with $\overline{P_{\infty}P_1} \ne \overline{P_{\infty}P_2}$. 
Then, $\ell$ is given by the intersection of planes spanned by $\overline{P_{\infty}P_1}$ and $R'$, $\overline{P_{\infty}P_2}$ and $R'$. 
This implies that $\ell$ is the tangent line at $P_{\infty}$. 
This is a contradiction.

We consider the case where $\ell \cap \{Z=0\}=\{Q\} \not\subset \mathcal{X}$. 
By Lemma \ref{hermitian}, ramification indices for each points of $\mathcal{X} \cap \{Z=0\}$ are $d$ or $d/(q+1)$. 
It follows from Lemma \ref{order} that $d=q^2-q+1$ or $d=q+1$. 
Assume that $d=q^2-q+1$ and $q>2$. 
Then, for each $P \in \mathcal{X} \cap \{Z=0\}$, the line $\overline{QP}$ intersects $\mathcal{X}$ at a unique point $P$. 
Then, $\pi_{\ell}$ is ramified at each point $P \in \mathcal{X} \cap \{Z=0\}$ and hence, the tangent line at $P$ intersects $\ell$. 
If $\ell \not\ni R'$, then the tangent lines are included in the plane spanned by $\ell$ and $R'$. 
This is a contradiction. 
Therefore, $\ell \ni R'$. 
We can assume that $Q$ is not $\mathbb{F}_{q^2}$-rational and $\overline{R'Q} \cap \mathcal{X}$ consists of $q^2-q+1$ points.
Then, $d=(q^3+1)-(q^2-q+1)=q^2-q+1$. 
This is a contradiction.  
Assume that $q=2$ and $d=q^2-q+1=3$. 
If there exist two lines in $\{Z=0\}$ containing $q+1$ points of $\mathcal{X}$, then $Q$ is $\mathbb{F}_{q^2}$-rational, since such lines are $\mathbb{F}_{q^2}$-rational. 
We can assume that $Q=R$. 
Then, $\pi_R(\ell) \in {\rm Sing}(\mathcal{X}') \setminus \{Z=0\}$. 
This is a contradiction. 
Therefore, there exist points $P_1$ and $P_2 \in \mathcal{X} \cap \{Z=0\}$ with $\overline{QP_1} \ne \overline{QP_2}$ such that the tangent lines at $P_1$ and $P_2$ intersect $\ell$.   
If $R' \not\in \ell$, then $P_1$ and $P_2$ are included in the plane spanned by $\ell$ and $R'$. 
This is a contradiction to that points $Q$, $P_1$ and $P_2$ are not collinear in the plane $\{Z=0\}$. 
Therefore, $\ell \ni R'$. 
We can assume that $Q$ is not $\mathbb{F}_{q^2}$-rational and $\overline{R'Q} \cap \mathcal{X}$ consists of $q^2-q+1$ points.
Then, $d=(q^3+1)-(q^2-q+1)=q^2-q+1$. 
This is a contradiction.

Assume that $d=q+1$. 
If $q=2$, then $d=q+1=q^2-q+1$. 
We can assume that $q>2$.  
It follows from Lemma \ref{order} that, for each $P \in \mathcal{X} \cap \{Z=0\}$, the line $\overline{QP}$ contains exactly $q+1$ points of $\mathcal{X}$. 
This implies that $\overline{QP}$ is $\mathbb{F}_{q^2}$-rational and hence, $Q$ is $\mathbb{F}_{q^2}$-rational. 
We can assume that $Q=R$. 
Then, $\pi_R(\ell) \in {\rm Sing}(\mathcal{X}') \setminus \{Z=0\}$. 
This is a contradiction.


\begin{thebibliography}{100} 
\bibitem{beelen-montanucci} P. Beelen and M. Montanucci, Weierstrass semigroups on the Giulietti--Korchm\'{a}ros curve, Finite Fields Appl. {\bf 52} (2018), 10--29. 
\bibitem{duyaguit-yoshihara} C. Duyaguit and H. Yoshihara, Galois lines for normal elliptic space curves, Algebra Colloq. {\bf 12} (2005), 205--212.
\bibitem{fukasawa} S. Fukasawa, A birational embedding of an algebraic curve into a projective plane with two Galois points, preprint, arXiv:1611.03953.
\bibitem{giulietti-korchmaros} M. Giulietti and G. Korchm\'{a}ros, A new family of maximal curves over a finite field, Math. Ann. {\bf 343} (2009), 229--245.  
\bibitem{gqz} M. Giulietti, L. Quoos and G. Zini, Maximal curves from subcovers of the GK-curve, J. Pure Appl. Algebra {\bf 220} (2016), 3372--3383. 
\bibitem{homma} M. Homma, Galois points for a Hermitian curve, Comm. Algebra {\bf 34} (2006), 4503--4511.
\bibitem{miura-yoshihara} K. Miura and H. Yoshihara, Field theory for function fields of plane quartic curves, J. Algebra {\bf 226} (2000), 283--294.
\bibitem{stichtenoth} H. Stichtenoth, Algebraic Function Fields and Codes, Universitext, Springer-Verlag, Berlin (1993).
\bibitem{suzuki} M. Suzuki, Group Theory I, Grundlehren der Mathematischen Wissenschaften {\bf 247}, Springer-Verlag, Berlin-New York (1982).
\bibitem{yoshihara} H. Yoshihara, Function field theory of plane curves by dual curves, J. Algebra {\bf 239} (2001), 340--355.
\bibitem{yoshihara2} H. Yoshihara, Galois lines for space curves, Algebra Colloq. {\bf 13} (2006), 455--469.  
\bibitem{yoshihara-fukasawa} H. Yoshihara and S. Fukasawa, List of problems, available at: \\
http://hyoshihara.web.fc2.com/openquestion.html
\end{thebibliography}
\end{document}